\documentclass[12pt]{article}
\usepackage{fullpage,graphicx,graphics,psfrag,amsmath,amsfonts,verbatim}
\usepackage[small,bf]{caption}
\usepackage{url,ar}
\usepackage{amssymb,enumitem}
\usepackage[siunitx]{circuitikz}

\newcommand{\BEAS}{\begin{eqnarray*}}
\newcommand{\EEAS}{\end{eqnarray*}}
\newcommand{\BEA}{\begin{eqnarray}}
\newcommand{\EEA}{\end{eqnarray}}
\newcommand{\BEQ}{\begin{equation}}
\newcommand{\EEQ}{\end{equation}}
\newcommand{\BIT}{\begin{itemize}}
\newcommand{\EIT}{\end{itemize}}
\newcommand{\BNUM}{\begin{enumerate}}
\newcommand{\ENUM}{\end{enumerate}}

\newcommand{\BA}{\begin{array}}
\newcommand{\EA}{\end{array}}

\newcommand{\eg}{{\it e.g.}}
\newcommand{\ie}{{\it i.e.}}

\newcommand{\reals}{{\mbox{\bf R}}}

\newcommand{\K}{\mathcal{K}}

\newtheorem{theorem}{Theorem}[section]

\newtheorem{claim}[theorem]{Claim}

\makeatletter
\long\def\@makecaption#1#2{
   \vskip 9pt
   \begin{small}
   \setbox\@tempboxa\hbox{{\bf #1:} #2}
   \ifdim \wd\@tempboxa > 5.5in
        \begin{center}
        \begin{minipage}[t]{5.5in}
        \addtolength{\baselineskip}{-0.95pt}
        {\bf #1:} #2 \par
        \addtolength{\baselineskip}{0.95pt}
        \end{minipage}
        \end{center}
   \else
    \hbox to\hsize{\hfil\box\@tempboxa\hfil}
   \fi
   \end{small}\par
}
\makeatother

\newcounter{oursection}

\newcounter{lecture}

\renewcommand{\K}{\mathcal{K}}
\newcommand{\F}{\mathcal{F}}

\newcounter{algorithmctr}[section]
\renewcommand{\thealgorithmctr}{\thesection.\arabic{algorithmctr}}
    {\refstepcounter{algorithmctr}\begin{list}{}{%
		\setlength{\rightmargin}{0\linewidth}%
		\setlength{\leftmargin}{.05\linewidth}}%
		\rmfamily\small
		\item[]{\setlength{\parskip}{0ex}\hrulefill\par%
		\nopagebreak{\bfseries\textsf{Algorithm \thealgorithmctr~}}}}%
	{{\setlength{\parskip}{-1ex}\nopagebreak\par\hrulefill} \end{list}}

\title{Disciplined Multi-Convex Programming}
\author{Xinyue Shen \and Steven Diamond \and Madeleine Udell
\and Yuantao Gu \and Stephen Boyd}

\begin{document}
\maketitle

\begin{abstract}
	A multi-convex optimization problem is one in which the variables 
	can be partitioned into sets over which the problem is convex 
	when the other variables are fixed.
	Multi-convex problems are generally solved approximately 
	using variations on alternating or cyclic minimization.  
	Multi-convex problems arise in many applications,
	such as nonnegative matrix factorization, generalized low rank models,
	and structured control synthesis, to name just a few.  
	In most applications to date the multi-convexity 
	is simple to verify by hand.
	In this paper we study the automatic detection
	and verification of multi-convexity using the ideas of 
	disciplined convex programming.
	We describe an implementation of our proposed method that
	detects and verifies multi-convexity, and then invokes
	one of the general solution methods.
\end{abstract}


\section{Introduction}

A multi-convex optimization problem is one in which the variables 
can be partitioned into sets over which the problem is convex 
when the other variables are fixed.
Multi-convex problems appear in domains such as
machine learning \cite{lee1999learning, udell2014generalized}, 
signal and information processing \cite{Lee00, Kim2008, Wen2012},
communication \cite{1254038}, 
and control \cite{751690, hassibi1999path, hours2014parametric, 7170815}.
Typical problems in these fields include nonnegative 
matrix factorization (NMF) and bilinear matrix inequality (BMI) problems.


In general multi-convex problems are hard to solve globally,
but several algorithms have been proposed as heuristic or local methods,
and are widely used in applications.
Most of these methods are variations on 
the block coordinate descent (BCD) method.
The idea of optimizing over a single block of variables while 
holding the remaining 
variables fixed in each iteration dates back to \cite{Warga,Powell1973}.
Convergence results were first discussed for strongly convex differentiable
objective function \cite{Warga}, and then under various assumptions
on the separability and regularity of the objective function 
\cite{Tseng1993,Tseng2001}.
In \cite{attouch2010}, a two-block BCD method with proximal operator 
is used to minimize a nonconvex objective function 
which satisfies the Kurdyka-Lojasiewicz inequality.
In \cite{razaviyayn2012unified} the authors propose an inexact BCD approach which 
updates variable blocks by minimizing a sequence of approximations 
of the objective function, which can either be nondifferentiable or nonconvex.
A recent work \cite{xu2013block} uses BCD to solve multi-convex 
problems, where the objective is a sum of a differentiable 
multi-convex function and several extended-valued convex functions.
In each step, updates with and without proximal operator 
and prox-linear operator are considered,
and convergence analysis is established under certain assumptions.
Gradient methods have also been proposed for multi-convex problems,
where the objective is differentiable in each block of variables,
and all variables are updated at once along their descent directions
and then projected into a convex feasible set in every iteration \cite{Chu04}.

The focus of this paper is not on solution methods,
but on a \emph{modeling framework} for expressing multi-convex
problems in a way that verifies the multi-convex structure, and 
can expose the structure to whatever solution algorithm is then used.
Modeling frameworks have been developed for convex problems,
\eg,  \texttt{CVX} \cite{cvx}, \texttt{YALMIP} \cite{Lofberg:04},
\texttt{CVXPY} \cite{cvxpy_paper}, and \texttt{Convex.jl} \cite{cvxjl}. 
These frameworks provide a uniform method for specifying convex
problems based on the idea of disciplined convex programming 
(DCP) \cite{GBY:06}.
This gives a simple method for verifying that a problem is convex, and
for automatically canonicalizing to a standard
generic form such as a cone program.
The goal of DCP (and these software frameworks) is not to detect or
determine feasibility of an arbitrary problem, but rather to
give a very simple set of rules that can be used to construct 
convex problems.

In this paper we extend the idea of DCP to multi-convex problems.
We propose a disciplined multi-convex programming (DMCP) rule set,
an extension of the DCP rule set.
Problem specifications that conform to the DMCP rule set can be verified as
convex in a group of variables, for any fixed values of the other variables,
using ideas that extend those in DCP.
We describe an efficient algorithm that can carry out the analysis of 
convexity of problem when an arbitrary group of variables is fixed at any
value.
As with DCP, the goal of DMCP is not to analyze multi-convexity of 
an arbitrary problem, but rather to give a simple set of rules which if followed
yields multi-convex problems.
In applications to date, such as NMF, verification of multi-convexity is simple, and 
can be done by hand or just simple observation.
With DMCP a far larger class of multi-convex problems can be constructed
in an organized way.

We describe a software implementation of the ideas developed in this paper,
called \texttt{DMCP}, a Python package that extends \texttt{CVXPY}.
It implements the DMCP verification and analysis methods, and then
heuristically solves a conforming problem via BCD type algorithms,
which we extend for general use to include slack variables to handle 
infeasibility.  A similar package, \texttt{MultiConvex}, has been
developed for the Julia package \texttt{Convex.jl}.
We illustrate the framework on a number of examples.

In \S\ref{s-multicvx} we carefully define multi-convexity of a function,
constraint, and problem. 
We review block coordinate descent methods, introducing new variants
with slack variables and generalized inequalities, in \S\ref{s-bcd}.
In \S\ref{s-dmcp} we describe the main ideas of DMCP and an efficient 
algorithm for verifying that a problem specification conforms to DMCP.
In \S\ref{s-impl} we describe our implementation of the package 
\texttt{DMCP}.
Finally, in \S\ref{s-examples} we describe a number of numerical examples.
Our goal there is not to show competitive results, in terms of 
solution quality or solve time, but rather to show the simplicity with which
the problem is specified, along with results that are at least 
comparable to those obtained with custom solvers for the specific problem.

\section{Multi-convex programming}\label{s-multicvx}
\subsection{Multi-convex function}

\paragraph{Fixing variables in a function.}
Consider a function $f:\reals^n\to\reals\cup\{\infty\}$,
and a partition of the variable $x\in \reals^n$ into blocks of variables
\[
x=(x_1,\ldots, x_N), \quad x_i\in\reals^{n_i},
\]
so $\sum_{i=1}^N n_i = n$.
Throughout this paper we will use subsets of indices to refer to 
sets of the variables.
Let $\F\subseteq\{1,\ldots,N\}$ denote an index set,
with complement $\F^c = \{1,\ldots,N\} \setminus \F$.
By fixing the variables with indices in $\F$
of the function $f$ at a given point $\hat{x}\in\reals^n$,
we obtain a function over the remaining variables, with indices in $\F^c$,
which we denote as $\tilde f = \mathrm{fix}(f, \hat{x}, \F) $.
For $i \in \F^c$, $x_i$ is a variable of the function $\tilde f$;
for $i \in \F$, $x_i=\hat x_i$.
Informally we refer to $\tilde f$ as 
`$f$, with the variables $x_i$ for $i\in \F$ fixed'.

As an example consider $f:\reals^4 \to \reals \cup \{\infty\}$ defined as
\begin{equation}\label{e-example}
f(x_1, x_2, x_3,x_4) = |x_1x_2 + x_3x_4|.
\end{equation}
With $\hat{x} = (1,2,1,3)$ and $\F=\{1,3\}$,
the fixed function $\tilde f = 
\mathrm{fix}(f, \hat{x}, \{1,3\})$ is given by
$\tilde f(x_2, x_4)  = |x_2+x_4|$.

\paragraph{Multi-convex and multi-affine functions.}
Given an index set $\F\subseteq\{1,\ldots,N\}$,
we say that a function $f:\reals^n\to\reals\cup\{\infty\}$
is \textit{convex} (or \textit{affine}) \textit{with set $\F$ fixed},
if for any $\hat{x}\in\reals^n$
the function $\mathrm{fix}(f,\hat{x},\F)$ is a convex (or affine) function.
(In this definition, we consider a so-called improper function
\cite{rockafellar}, which has the value $\infty$ everywhere,
as convex or affine.)
For example, the function $f$ defined in \eqref{e-example} 
is convex with the variables $x_1$ and $x_3$ fixed (\ie,
with index set $\F=\{1,3\}$).
A function is convex if and only if it is convex 
with $\F = \emptyset$ fixed,
\ie, with none of its variables fixed.

We say the function $f$ is \textit{multi-convex} 
(or \textit{multi-affine)}, if there are index sets $\F_1,\ldots,\F_K$,
such that for every $k$ the function $f$ is convex (or affine) with $\F_k$ fixed, 
and $\cap_{k=1}^K \F_k = \emptyset$.
The requirement that $\cap_{k=1}^K \F_k = \emptyset$ means that for 
every variable $x_i$ there is some $\F_k$ with $i \not\in \F_k$.
In particular, $\tilde{f}= \mathrm{fix}(f,\hat{x},\{1,\ldots,N\}\setminus \{i\})$
is convex in $x_i$.
For $K=2$, we say that the function is \emph{bi-convex} 
(or \emph{bi-affine}).

As an example, the function $f$ in \eqref{e-example} is 
convex with $\{1,3\}$ fixed and $\{2,4\}$ fixed, so it is multi-convex.
The choice of $K$, and the index sets, is not unique. 
The function $f$ is also convex
with $\{2,3,4\}$ fixed, $\{1,3,4\}$ fixed, $\{1,2,4\}$ fixed, and $\{1,2,3\}$ fixed.

\paragraph{Minimal fixed sets.}
For a function $f$, we can consider the set of all index sets $\F$ for which
$\mathrm{fix}(f, \hat x, \F)$ is convex for all $\hat x$;
among these we are interested in
the \emph{minimal fixed sets} that render a function convex.
A minimal fixed set is a set of variables that when fixed make the 
function convex; but if any variable is removed from the set, the function 
is not convex.
A function is multi-convex if and only if the intersection of these
minimal fixed index sets is empty.

\subsection{Multi-convex problem}
We now extend the idea of multi-convexity to the optimization problem
\begin{equation}\label{main}
\begin{array}{ll}
\mbox{minimize} &  f_0(x) \\
\mbox{subject to} & f_i (x) \leq 0, \quad i = 1,\ldots, m \\
& g_i (x)  = 0, \quad i = 1,\ldots, p,
\end{array}
\end{equation}
with variable $x\in\reals^n$ partitioned into blocks as
$x = (x_1,\ldots,x_N)$,
and functions $f_i:\reals^{n}\to\reals\cup\{\infty\}$ for $i=0,\ldots,m$
and $g_i:\reals^n\to\reals$ for $i=1,\ldots, p$ are proper.

Given an index set $\F\subseteq\{1,\ldots,N\}$,
problem~\eqref{main} is \textit{convex with set $\F$ fixed},
if for any $\hat{x}\in\reals^n$ the problem
\begin{equation}\label{fix_p}
\begin{array}{ll}
\mbox{minimize} &  \mathrm{fix}(f_0, \hat{x}, \F) \\
\mbox{subject to} &\mathrm{fix} (f_i, \hat{x}, \F) \leq 0, \quad i = 1,\ldots, m \\
& \mathrm{fix} (g_i, \hat{x}, \F)   = 0, \quad i = 1,\ldots, p,
\end{array}
\end{equation}
is convex.
In other words, problem~\eqref{main} is convex with $\F$ fixed,
if and only if functions $f_i$ for $i=0,\ldots, m$ are convex with $\F$ fixed,
and functions $g_i$ for $i=1,\ldots,p$ are affine with $\F$ fixed.

We say the problem~\eqref{main} is \textit{multi-convex}, if there are sets
$\F_1,\ldots, \F_K$,
such that for every $k$ problem~\eqref{main} is convex with set $\F_k$ fixed, 
and $\cap_{k=1}^K \F_k = \emptyset$.
A convex problem is multi-convex with $K=0$ (\ie, $\F = \emptyset$).
A \emph{bi-convex problem} is multi-convex with $K=2$.

As an example the following problem is multi-convex:
\begin{equation}\label{example}
\begin{array}{ll}
\mbox{minimize} & |x_1x_2 + x_3x_4| \\
\mbox{subject to} & x_1+x_2+x_3+x_4 = 1,
\end{array}
\end{equation}
with variable $x\in \reals^4$.
This is readily verfied with 
$\mathcal F_1 = \{1,3\}$ and $\mathcal F_2 = \{2,4\}$.
For a given problem we can consider the minimal variable index
sets which make the problem convex.  If the problem is convex, 
$\F = \emptyset$ is the unique minimal set.

\section{Block coordinate descent and variations}\label{s-bcd}

In this section we review, and extend, some generic methods for approximately
solving the multi-convex problem \eqref{main}, using BCD-type methods.

\subsection{Block coordinate minimization with slack variables}
Assume that sets $\F_k$, $k=1,\ldots,K$ are index sets
for which the problem~\eqref{main} with $\F_k$ fixed is convex, 
with $\cap_{k=1}^K \mathcal F_k =\emptyset$.
These could be the set of all minimal index sets, but any other set of index sets
that verify multi-convexity could be used.

The basic form of the proposed method is iterative.  In each iteration,
we fix the variables in one set $\F_k$ and solve the following subproblem,
\begin{equation}\label{slack}
\begin{array}{ll}
\mbox{minimize} &  \mathrm{fix}( f_0, \hat{x}, \F_k) 
+\mu \sum_{i=1}^{m} s_i + \mu\sum_{i=1}^p |s_{i+m}|\\
\mbox{subject to} &  \mathrm{fix}( f_i, \hat{x}, \F_k) \leq s_i, \quad s_i\geq0,\quad i = 1,\ldots, m \\
&  \mathrm{fix}( g_i, \hat{x}, \F_k)  = s_{i+m}, \quad i = 1,\ldots, p, \\
\end{array}
\end{equation}
where $s_i$ for $i=1,\ldots,m+p$ and $x_i$ for $i\in\F_k^c$ are the variables,
and $\mu>0$ is a parameter.
Here the constant $\hat{x}$ inherits the value of $x$ from the least iteration.
This subproblem solved in each iteration is convex. 
The slack variables $s_i$ for $i=1,\ldots,m+p$ ensure that the subproblem 
\eqref{slack} cannot be infeasible.
The added terms in the objective are a so-called exact penalty 
\cite{nocedal2006numerical},
meaning that when some technical conditions hold, and $\mu$ is large enough, 
the solution satisfies $s_i=0$, when the subproblem without the slack 
variables is feasible.

Many schemes can be used to choose $k$ in each iteration, and to update the 
slack parameter $\mu$.
For example, we can cyclically choose $k$, or randomly choose $k$, 
or optimize over $k=1, \ldots, K$ in rounds of $K$ steps,
in an order chosen by a random permutation in each round.
Updating $\mu$ is typically done by increasing it by a factor $\rho>1$ after each
iteration, or after each round of $K$ iterations.
One variation on the algorithm sets the slack variables to zero (\ie, removes
them) once a feasible point is obtained (\ie, a point is obtained with
$s_i=0$).
The algorithm is typically initialized with values specific to the particular 
application, or generic values.

This algorithm differs from the basic BCD algorithm in the addition
of the slack variables, and in the feature that a variable can appear in more
than one set $\F_k^c$, meaning that a variable can be updated
in multiple iterations per round of $K$ iterations.
For example, if a problem has variables $x_1,x_2,x_3$ and is convex in
$(x_1,x_2)$ and $(x_2,x_3)$, our method will update $x_2$ in each step.

In the general case, very little can be said about the convergence of this method.
One obvious observation is that, if $\mu$ is held fixed, the objective
is nonincreasing and so convergences.
See the references cited above for some convergence results for related algorithms,
for special cases with strong assumptions such as strict convexity (when the variables
are fixed) or differentiability.
As a practical matter, similar algorithms have been found to be robust, and 
very useful in practice, despite a lack of strong theory
establishing convergence in the general case.

\subsection{Block coordinate proximal iteration}
A variation of subproblem \eqref{slack} is adds a proximal term 
\cite{parikh2014proximal},
which renders the subproblems strongly convex:
\begin{equation}\label{proximal}
\begin{array}{ll}
\mbox{minimize} &   \mathrm{fix}( f_0, \hat{x}, \F_k)
+\mu \sum_{i=1}^{m} s_i + \mu \sum_{i=1}^p |s_{i+m}| 
+\frac{1}{2\lambda} \sum_{i\in\F^c_k} \|x_i-\hat{x}_i\|_2^2\\
\mbox{subject to} 
& \mathrm{fix}( f_i, \hat{x}, \F_k) \leq s_i, \quad s_i\geq0, \quad i = 1,\ldots, m \\
&  \mathrm{fix}( g_i, \hat{x}, \F_k) = s_{i+m}, \quad i = 1,\ldots, p,
\end{array}
\end{equation}
where $x_i$ for $i\in\F^c_k$ and $s_i$ for $i=1,\ldots, m+p$ are variables,
$\lambda >0$ is the proximal parameter.
The proximal term penalizes large changes in the variables being 
optimized, \ie, it introduces damping into the algorithm.  In some cases
it has been observed to yield better final points, \ie, points
with smaller objective value, than those obtained without proximal
regularization.

Yet another variation uses linearized proximal steps, when $f$ is differentiable
in the variables $x_i$ for $i \in \mathcal F_k^c$.  The subproblem
solved in this case is
\begin{equation}\label{linear}
\begin{array}{ll}
\mbox{minimize} & 
\mu \sum_{i=1}^{m} s_i + \mu\sum_{i=1}^p |s_{i+m}| 
+ \sum_{i\in\F^c_k} \left(  \frac{1}{2\lambda} \|x_i-\hat{x}_i\|_2^2 
+(x_i-\hat{x}_i)^T \nabla f(\hat{x}_i) \right) \\
\mbox{subject to} 
&  \mathrm{fix}( f_i, \hat{x}, \F_k)\leq s_i, \quad s_i\geq0, \quad i = 1,\ldots, m \\
& \mathrm{fix}( g_i, \hat{x}, \F_k)  = s_{i+m}, \quad i = 1,\ldots, p,
\end{array}
\end{equation}
where $x_i$ for $i\in\F^c_k$ and $s_i$ for $i=1,\ldots, m+p$ are variables,
and $\nabla f(\hat{x}_i)$ is the partial gradient of $f$ with respect to $x_i$ at
the point $\hat{x}$.
The objective is equivalent to the minimization of 
\[
\mu \sum_{i=1}^{m} s_i + \mu\sum_{i=1}^p |s_{i+m}| 
+\sum_{i\in\F^c_k} \frac{1}{2\lambda} \|x_i-\hat{x}_i + 
\lambda \nabla f(\hat{x}_i)\|_2^2, 
\]
which is the objective of a proximal gradient method.

\subsection{Generalized inequality constraints}
One useful extension is to generalize problem~\eqref{main} by allowing
generalized inequality constraints.
Suppose the functions $f_0$ and $g_i$ for $i=1,\ldots, p$ are the same as
in problem~\eqref{main}, 
but $f_i : \reals^n\to\reals^{d_i}\cup\{\infty\}$.
Consider the following program with generalized inequalities,
\begin{equation}
\begin{array}{ll}
\mbox{minimize} & f_0(x) \\
\mbox{subject to} & f_i (x) \preceq_{\K_i} 0,\quad i=1,\ldots, m \\
& g_i(x) = 0,\quad i = 1,\ldots,p,
\end{array}
\end{equation}
where $x = (x_1,\ldots, x_N)\in\reals^n$ is the variable, 
and the generalized inequality constraints
are with respect to proper cones $\K_i\subseteq\reals^{d_i}$, 
$i=1, \ldots, m$.
The definitions of multi-convex program and minimal index set 
can be directly extended.
Slack variables are added in the following way:
\begin{equation}\label{ext}
\begin{array}{ll}
\mbox{minimize} &  \mathrm{fix}( f_0, \hat{x}, \F_k) 
+\mu \sum_{i=1}^{m} s_i + \mu\sum_{i=1}^p |s_{i+m}|\\
\mbox{subject to} 
&  \mathrm{fix}( f_i, \hat{x}, \F_k) \preceq_{\K_i} s_i e_i, 
 \quad s_i\geq0,\quad i = 1,\ldots, m \\
&  \mathrm{fix}( g_i, \hat{x}, \F_k)  = s_{i+m}, \quad  i = 1,\ldots, p,
\end{array}
\end{equation}
where $e_i$ is a given positive element in cone $\K_i$
for $i=1,\ldots,m$.

\section{Disciplined multi-convex programming}\label{s-dmcp}
\subsection{Disciplined convex programming}
Disciplined convex programming (DCP) is a methodology introduced by Grant et al.
\cite{GBY:06}
that imposes a set of conventions that must be followed when
constructing (or specifying or defining) convex programs.
Conforming problems are called \emph{disciplined convex programs}.
A disciplined convex program can be transformed into an
equivalent cone program by replacing each function 
with its graph implementation \cite{GB:08}.
The convex optimization modeling systems
\texttt{YALMIP} \cite{Lofberg:04}, \texttt{CVX} \cite{cvx},
\texttt{CVXPY} \cite{cvxpy_paper},
and \texttt{Convex.jl} \cite{cvxjl} use DCP
to verify the convexity of a problem
and automatically convert convex programs into cone programs,
which can then be solved using generic solvers.

The conventions of DCP restrict the set of functions that can appear in
a problem and the way functions can be composed.
Every function in a disciplined convex program must be formed 
as an expression involving constants or parameters, variables,
and a dictionary of atomic functions.
The dictionary consists of functions with known curvature and monotonicity,
and a graph implementation, or representation as partial
optimization over a cone program \cite{BoV:04,NesNem:92}.
Every composition of functions $f(g_1(x),\ldots,g_p(x))$,
where $f : \reals^p \to \reals$ is convex and
$g_1,\ldots,g_p : \reals^n \to \reals$,
must satisfy the following composition rule,
which ensures the composition is convex.
Let $\tilde{f} : \reals^p \to \reals \cup \{\infty\}$
be the extended-value extension of $f$ \cite[Chap.~3]{BoV:04}.
One of the following conditions must hold for each $i=1,\ldots,p$:
\begin{itemize}
	\item $g_i$ is convex and $\tilde{f}$ is nondecreasing in argument $i$.
	\item $g_i$ is concave and $\tilde{f}$ is nonincreasing in argument $i$.
	\item $g_i$ is affine.
\end{itemize}
The composition rule for concave functions is analogous.

Signed DCP is an extension of DCP that keeps track of the 
signs of functions and expressions, using simple sign arithmetic.
The monotonicity of functions in the atom library can then depend
on the sign of their arguments.
As a simple example, 
consider the expression $y=(\exp x)^2$, where $x$ is a variable.
The subexpression $\exp x$ is convex and (in signed DCP analysis) 
nonnegative.  With DCP analysis, $y$ cannot be verified as
convex, since the square function is not nondecreasing. With signed
DCP analysis, the square function is known to be nondecreasing
for nonnegative arguments, which matches this case, so $y$ is 
verified as convex using signed DCP analysis.

Convexity verification of an expression formed from variables 
and constants (or parameters) in (signed) DCP first analyzes the signs
of all subexpressions.  Then it checks that the composition rule
above holds for every subexpression, possibly relying on the known
signs of subexpressions.
If everything checks out the expression is verified to be constant,
affine, convex, concave, or unknown (when the DCP rules do not hold).
We make an observation that is critical for our work here:
The DCP analysis does not use the values of any constants or parameters
in the expression.  The number 4.57 is simply treated as positive;
if a parameter has been declared as positive, then it is treated
as positive.
It follows immediately that DCP analysis has verified not just that
the specific expression is convex, but that it is convex for
\emph{any} other values of the constants and parameters, with the 
same signs as the given ones, if the sign matters.

The following code snippet gives an example in \texttt{CVXPY}:
\begin{quote}
	\begin{verbatim}
	x = Variable(n)
	mu = Parameter(sign = 'positive')
	expr = sum_squares(x) + mu*norm(x,1)
	\end{verbatim}
\end{quote}
In the first line we declare (or construct) a variable, 
and in the second line construct a paramater, \ie, a constant
that is unknown, but positive.
The curvature of the expression \texttt{expr} is verified to be convex,
even though the value of parameter \texttt{mu} is unknown;
DCP analysis uses only the fact that whatever the value of $\texttt{mu}$
is, it must be positive.

\subsection{Disciplined multi-convex programming}

To determine that problem~\eqref{main} 
is multi-convex requires us to verify that functions
$\mathrm{fix}(f_i,\hat{x},\F)$ are convex,
and if $\mathrm{fix}(g_i,\hat{x},\F)$ are affine, for all $\hat{x}\in\reals^n$.
We can use the idea of DCP, specifically with signed parameters, 
to carry this out,
which gives us a practical method for multi-convexity verification.

\paragraph{Multi-convex atoms.}
We start by generalizing the library of DCP atom functions to include 
multi-convex atomic functions.
A function is a multi-convex atom if it has $N$ arguments $N>1$,
and it reduces to a DCP atomic function, when and only when 
all but the $i$th arguments are constant for each $i=1,\ldots,N$.
For example, the product of $N$ variables is a multi-convex atom
that extends the DCP atom of multiplication between one variable and constants.

Given a description of problem~\eqref{main} under a library of DCP and 
multi-convex atomic functions,
we say that it is \textit{disciplined convex programming 
	with set $\F\subseteq\{1,\ldots,N\}$ fixed},
if the corresponding problem~\eqref{fix_p} for any $\hat{x}\in\reals^n$
conforms to the DCP rules with respect to the DCP atomic function set.
When there is no confusion,
we simply say that problem~\eqref{main} is DCP with $\F$ fixed.

To verify if a problem is DCP with $\F$ fixed,
a method first fixes the problem by replacing variables in $\F$ 
with parameters of the same signs and dimensions.
Then it verifies DCP of the fixed problem with parameters 
according to the DCP ruleset.
The parameter is the correct model for fixed variables,
in that the DCP rules ensure that the verified curvature holds for
any value of the parameter.

\paragraph{Disciplined multi-convex program.}
Given a description of problem~\eqref{main} under a library of DCP and 
multi-convex atomic functions,
it is \textit{disciplined multi-convex programming} (DMCP),  
if there are sets $\F_1,\ldots, \F_K$ such that
problem~\eqref{main} with every $\F_k$ fixed is DCP,
$\cap_{k=1}^K \F_k = \emptyset$.
We simply say that problem~\eqref{main} is DMCP if there is no confusion.
A problem that is DMCP is guaranteed to be multi-convex, just as
a problem that is DCP is guaranteed to be convex.
Morever, when a BCD method is applied to a DMCP problem, each iteration
involves the solution of a DCP problem.

\paragraph{DMCP verification.}
A direct way of DMCP verification is to check if the problem is DCP with
$\{i\}^c$ fixed for every $i=1,\ldots,N$.
Expressions in DMCP inherit the tree structure from DCP,
so such verification can be done in $O(MN)$ time,
where $M$ is number of nodes in problem expression trees, 
and $N$ is the number of distinct variables.
To see why DMCP can be verified in this simple way, 
we have the following claim.
The claim implies that every DMCP problem is DCP when all but one variable
is fixed.

\begin{claim}\label{claim:DCP_exp}
	For a problem consisting only of DCP atoms (or multi-convex atoms
	with all but one arguments constant),
	it is DCP with $\F$ fixed, 
	if and only if it is DCP with $\{i\}^c$ fixed for all $i\in\F^c$.
\end{claim}

To see why Claim~\ref{claim:DCP_exp} is correct, we first prove the direction 
that if a problem consisting only of DCP atoms is DCP with $\{i\}^c$ fixed 
for all $i\in\F^c$, then it is DCP with $\F$ fixed.
The proof begins with two observations for functions consisting only of DCP
atoms.
\begin{itemize}
	\item The DCP curvature types have a hierarchy:
	\[
	\mathrm{unknown}\to\left\{
	\begin{array}{ll}
	&\mathrm{convex} \\
	&\mathrm{concave}
	\end{array}\;\;
	\right\} \to \mathrm{affine} \to \mathrm{constant}. 
	\]
	Unknown is the base type, then it splits into convex and concave. 
	Affine is a subtype of both convex and concave. 
	Then constant is a subtype of everything. 
	There is a similar hierarchy for sign information.
	The DCP type system is monotone in the curvature and sign hierarchies, 
	meaning if the curvature or sign of an argument of a function is changed
	to be more specific, the type of the function will become more specific
	or stay the same.
	Fixing variables of a function 
	makes the curvatures of some arguments more specific, 
	while keeps all signs the same,
	so the curvature of the function can only get more specific or stay the same. 
	\item If a function is affine in $x$ and $y$ separately, 
	then it is affine in $(x,y)$. 
	This is true because no multiplication of variables is allowed in DCP,
	and the function can only be in the form of $Ax+By+c$.
\end{itemize}
Now suppose that a problem consisting of DCP atoms
with only two variables $x$ and $y$ is not DCP,
but it is DCP with $x$ fixed and with $y$ fixed.
Then there must be some function that has a wrong curvature type in $(x,y)$
but whose arguments all have known curvatures. 
Since the function has the right curvature type with $x$ fixed and $y$ fixed, 
there must be an argument that is convex or concave (not affine) in $(x,y)$, 
but has a different curvature in $x$ and $y$.
According to the first observation, the argument must be affine in $x$ and $y$. 
By the second observation, the argument is affine in $(x,y)$, 
which is a contradiction.
For the same reason, cases with more than two variables have the same conclusion.

For the other direction of Claim~\ref{claim:DCP_exp},
we again use the observation that if a problem is DCP with $\F$ fixed,
then fixing additional variables only makes function curvatures
more specific.
The problem must then be DCP with $\{i\}^c$ fixed
for all $i \in \F^c$, since $\F \subseteq \{i\}^c$.

%


\subsection{Efficient search for minimal sets}\label{s-min-sets}
A problem may have multiple collections of index sets for which it is DMCP.
We propose several generic and efficient ways of choosing which collection to use
when applying a BCD method to the problem.
The simplest option is to always choose the collection $\{1\}^c,\ldots,\{N\}^c$,
in which case BCD optimizes over one variable at a time.

A more sophisticated approach is to reduce the collection $\{1\}^c,\ldots,\{N\}^c$ to
minimal sets,
which allows BCD to optimize over multiple variables each iteration.
We find minimal sets by first determining which variables can be optimized
together.
Concretely, we construct a conflict graph $(\mathcal{V}, \mathcal{E})$,
where $\mathcal{V}$ is the set of all variables,
and $i\sim j\in\mathcal{E}$ if and only if variables $i$ and $j$
appear in two different child trees of a multi-convex atom in the problem
expression tree,
which means the variables cannot be optimized together.

Constructing the conflict graph takes $O(N^2M)$ time.
We simply do a depth-first traversal of the problem expression tree.
At each leaf node, we initialize a linked list with the leaf variable.
At each parent node, we join the linked lists of its children.
At each multi-convex atom node, we also remove duplicates from each child's
linked list and then iterate over the lists,
adding an edge for every two variables appearing in different lists.
The edges added at a given multi-convex node are all unique
because a duplicate edge would mean the same variable appeared in two different
child trees, which is not possible in a DMCP problem.
Hence, iterating over the lists of variables takes at most $N^2$ operations.

Given the conflict graph, for $i=1,\ldots,N$ we find a maximal independent
set $\F_i$ containing $i$ using a standard fast algorithm and replace $\{i\}^c$
with $\F_i^c$.
The final collection is all index sets $\F_i^c$ that are not supersets
of another index set $\F_j^c$.
More generally, we can choose any collection of index sets $\F_1,\ldots,\F_K$ such
that $\F_1^c,\ldots,\F_K^c$ are independent sets in the conflict graph.

\section{Implementation}\label{s-impl}
The methods of DMCP verification, searching for minimal sets to fix, 
and cyclic optimization with
minimal sets fixed are implemented as an extension of \texttt{CVXPY} 
in a package \texttt{DMCP} that can be accessed at
\verb|https://github.com/cvxgrp/dmcp|.
A Julia package with similar functionality, \texttt{MultiConvex.jl}, 
can be found at
\verb|https://github.com/madeleineudell/MultiConvex.jl|,
but we focus here on the Python package \texttt{DMCP}.


\subsection{Some useful functions}
\paragraph{Multi-convex atomic functions.}
In order to allow multi-convex functions,
we extend the atomic function set of \texttt{CVXPY}.
The following atoms are allowed to have non-constant expressions
in both arguments, 
while in base \texttt{CVXPY} one of the arguments must be constant.
\BIT
\item 
multiplication: \verb|expression1 * expression2|
\item 
elementwise multiplication: \verb|mul_elemwise(expression1, expression2)|
\item 
convolution: \verb|conv(expression1, expression2)|
\EIT

\paragraph{Find minimal sets.}
Given a problem, the function
\begin{quote}
	\begin{verbatim}
	find_minimal_sets(problem)
	\end{verbatim}
\end{quote}
runs the algorithm discussed in \S~\ref{s-min-sets} 
and returns a list of minimal sets of indices of variables.
The indices are with respect to the list \verb|problem.variables()|,
namely, the variable corresponding to index $0$ is
\verb|problem.variables()[0]|.

\paragraph{DMCP verification.}
Given a problem, the function
\begin{quote}
	\begin{verbatim}
	is_dmcp(problem)
	\end{verbatim}
\end{quote}
returns a boolean indicating if it is a DMCP problem.

\paragraph{Fix variables.}
The function
\begin{quote}
	\begin{verbatim}
	fix(expression, fix_vars)
	\end{verbatim}
\end{quote}
returns a new expression with the variables in the list \verb|fix_vars|
replaced with parameters of the same signs and values.
If \verb|expression| is replaced with a \verb|CVXPY| problem,
 then a fixed problem is returned by fixing every expression in its 
 cost function and both sides of inequalities or equalities in constraints.

\paragraph{Random initialization.}
It is suggested that users provide an initial point $x^0$ for the method such that
functions $\mathrm{fix}(f_i,x^0,\F_1)$ are proper for $i=0,\ldots,m$,
where $\F_1$ is the first minimal set given by \verb|find_minimal_sets|.
If not, the function \verb|rand_initial(problem)|
will be called to generate random values from the uniform distribution
over the interval $[0,1)$ ($(-1,0]$) for variables with non-negative (non-positive) sign,
and from the standard normal distribution for variables with no sign.
There is no guarantee
that such a simple random initialization can always work for any problem.

%
%

\subsection{Options on update and algorithm parameters}
The solving method is to cyclically fix every minimal set
found by \verb|find_minimal_sets| and update the variables.
Three ways of updating variables are implemented.
The default one can be called by
\verb|problem.solve(method = 'bcd', update = 'proximal')|,
which is to solve the subproblem with proximal operators,
\ie, problem~\eqref{proximal}. 
To update by minimizing the subproblem without proximal operators,
\ie, problem~\eqref{slack}, the solve method is called with
\verb|update = 'minimize'|.
To use the prox-linear operator in updates, \ie, problem~\eqref{linear},
the solve method should be called with
\verb|update = 'prox_linear'|.

The parameter $\mu$ is update in every cycle by
$\mu_{t+1} = \mathrm{min} (\rho\mu_t, \mu_\mathrm{max})$.
The algorithm parameters are $\rho,\mu_0, \mu_\mathrm{max}, \lambda$,
and the maximum number of iterations.
They can be set by passing values of the parameters
 \verb|rho|, \verb|mu_0|, \verb|mu_max|,
\verb|lambd|, and \verb|max_iter|, respectively, to the \verb|solve| method. 

\section{Numerical examples}\label{s-examples}
\subsection{One basic example}
\paragraph{Problem description.}
The first example is problem~\eqref{example},
which has appeared throughout this paper to explain definitions.
\paragraph{DMCP specification.}
The code written in \verb|DMCP| for this example is as follows.
\begin{quote}
	\begin{verbatim}
	prob = Problem(Minimize(abs(x_1*x_2+x_3*x_4)), [x_1+x_2+x_3+x_4 == 1])
	\end{verbatim}
\end{quote}
To find all minimal sets, the following line is typed in
\begin{quote}
	\begin{verbatim}
	find_minimal_sets(prob)
	\end{verbatim}
\end{quote}
and the output is \texttt{[[2, 1], [3, 1], [2, 0], [3, 0]]}. 
Note that index $i$ corresponds to variable \verb|prob.variables()[i]|
for $i=0,\ldots, 3$.
To verify if it is DMCP, the function
\begin{quote}
	\begin{verbatim}
	is_dmcp(prob)
	\end{verbatim}
\end{quote}
returns \verb|True|. 

\paragraph{Numerical result.}

Random initial values are set for all variables. The solve method with
default setting  finds a feasible point
with objective value $0$, which solves the problem globally.

\subsection{Fractional optimization}

\paragraph{Problem description.}
In this example we evaluate DMCP on some fractional optimization problems.
Consider the following problem
\begin{equation}\label{frac}
\begin{array}{ll}
\mbox{minimize} & p(x)/q(x) \\
\mbox{subject to} & x\in\mathcal{X},
\end{array}
\end{equation}
where $x\in\reals^n$ is the variable, $\mathcal{X}$ is a convex set, 
$p$ is a convex function, and $q$ is concave. The objective function
is set to $+\infty$ unless $p(x) \geq 0, q(x)>0$.
Such a problem is quasi-convex, and can be globally solved \cite[\S4.2.5]{BoV:04}, 
or even have analytical solutions, 
so the aim here is just to evaluate the effectiveness of the method.

There are several ways of specifying problem~\eqref{frac} as DMCP.
One way is via the following problem.
\begin{equation}\label{frac_t1}
\begin{array}{ll}
\mbox{minimize} & p(x)/q(y) \\
\mbox{subject to} & x\in\mathcal{X}, \quad x = y,
\end{array}
\end{equation}
where $x$ and $y$ are variables.
Another way is via the following.
\begin{equation}\label{frac_t2}
\begin{array}{ll}
\mbox{minimize} & \alpha \\
\mbox{subject to} & x\in\mathcal{X}, \quad p(x) \leq \alpha q(x),
\end{array}
\end{equation}
where $\alpha\in\reals_+$ and $x$ are variables.
Both of them are biconvex.

\paragraph{DMCP specification.}
Suppose that $\mathcal{X} = \reals^n$.
The code for the formulation in \eqref{frac_t1} is as the following.
\begin{quote}
	\begin{verbatim}
	x = Variable(n)
	y = Variable(n)
	# specify p and q here
	prob = Problem(Minimize(inv_pos(q)*p), [x == y])
	\end{verbatim}
\end{quote}
Expressions \texttt{p} and \texttt{q} are to be specified.
The code for problem~\eqref{frac_t2} is as follows.
\begin{quote}
	\begin{verbatim}
	alpha = Variable(1, sign = 'Positive')
	x = Variable(n)
	# specify p and q here
	prob = Problem(Minimize(alpha), [p <= q*alpha])	
	\end{verbatim}
\end{quote}

\paragraph{Numerical result.}
Take an example of $p(x) = x^2+1$ and $q(y) = \sqrt{y+0.5}$.
The code for specifying $p$ and $q$ is as the follows.
\begin{quote}
	\begin{verbatim}
p = square(x) + 1
q = sqrt(y+0.5)
	\end{verbatim}
\end{quote}
The global optimal value of the objective function is approximately $1.217$.
With random initial points, DMCP finds the global optimum for problem~\eqref{frac_t1}
 and \eqref{frac_t2}.

\subsection{Linear transceiver design}
\paragraph{Problem description.}
Suppose that a signal $x\in\reals^n$ passes through
a linear pre-coder $A\in\reals^{n\times n}$,
and is transmitted as $Ax$.
Denote the channel matrix as $C\in\reals^{m \times n}$
 and the additive noise as $e\in\reals^m$,
then the received signal is $y = CAx+e$.
The received signal passing through an equalizer $B\in\reals^{n \times m}$
is decoded as $By$.
The problem of determining $A$ and $B$ is called transceiver deign \cite{1254038}.

In this example, we assume that the signal $x$ is binary
and follows IID Bernoulli distribution, and the noise
$e\sim\mathcal{N}(0,\sigma_\mathrm{e}^2 I)$.
Given the channel matrix $C$, the aim is to design $A$ and $B$
such that the mean squared error $\mathbb{E}\|x-By\|_2^2$, 
where the mean is taken over $x$ and $e$, is minimized,
and the transmission power is constrained.

An optimization problem is formulated as the following
\[
\begin{array}{ll}
\mbox{minimize} & \frac{1}{2}\|BCA-I\|_F^2 +\sigma_\mathrm{e}^2 \|B\|_F^2 \\
\mbox{subject to} & \|A\|_F \leq p,
\end{array}
\]
where $B$ and $A$ are the variables. The problem is biconvex.

\paragraph{DMCP specification.}
The code can be written as the following.
\begin{quote}
\begin{verbatim}
A = Variable(n,n)
B = Variable(n,m)
sigma_e = Parameter(1)
cost = square(norm(B*C*A-I,'fro'))/2+square(sigma_e)*square(norm(B,'fro'))
prob = Problem(Minmize(cost), [norm(A, 'fro') <= p])
prob.solve(method = 'bcd')
\end{verbatim}
\end{quote}

\begin{figure}
	\centering
	\includegraphics[width=0.6\textwidth]{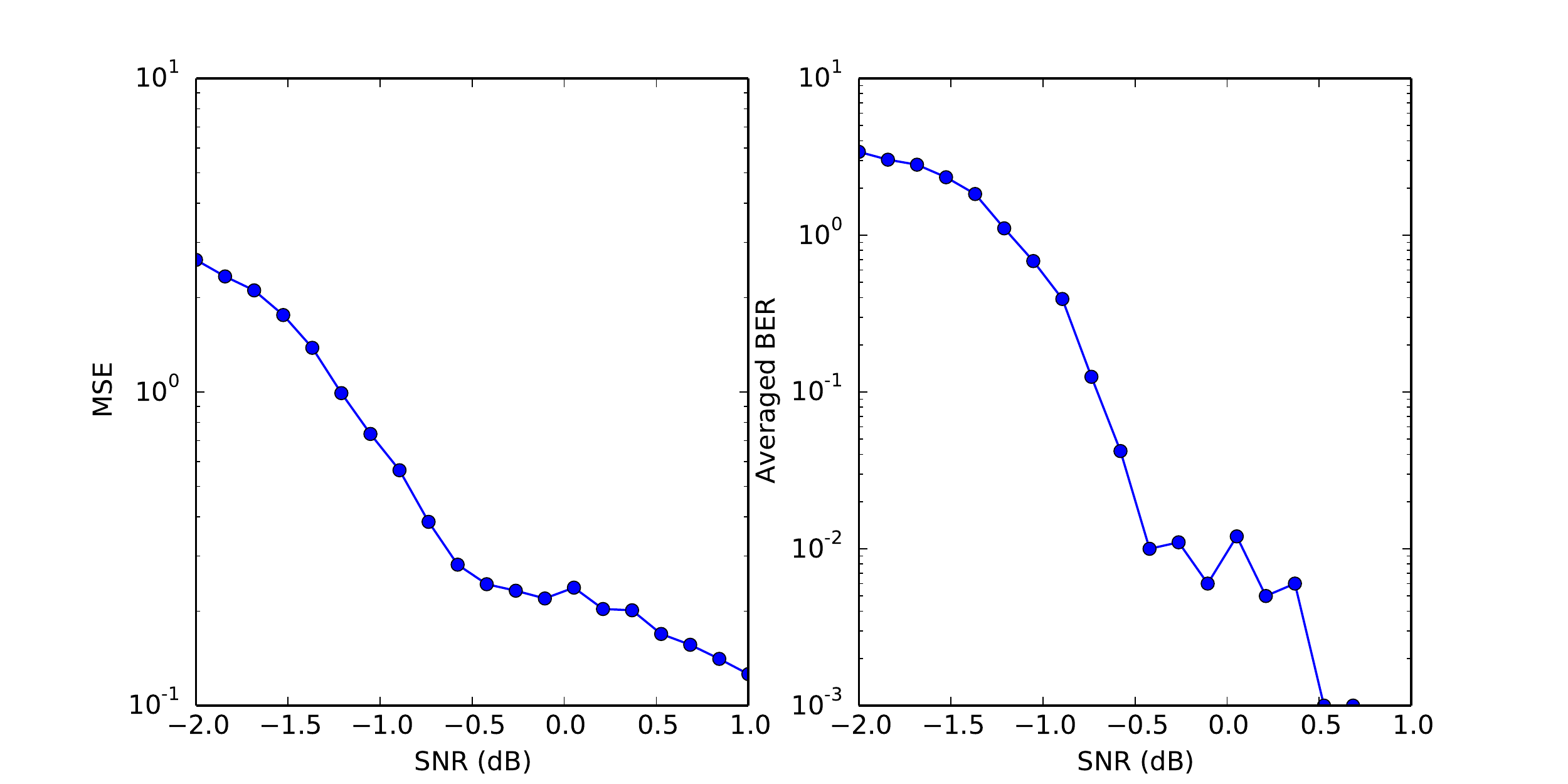}
	\caption{Linear transceiver design.}
	\label{fig:precoder_eq}
\end{figure}

\paragraph{Numerical result.}
In an experiment, $n = 10$, $m=15$, $p=10$, 
and the channel matrix $C$ is a random matrix with IID
normal distribution.
The signal to noise ratio varies, and for each value of $\sigma_\mathrm{e}$,
we try to solve the problem to get a design of $A$ and $B$.
Each design is tested by $1000$ trials with random signal $x$ and noise $e$
generated from the same distributions as the ones in the design.
The method is run without proximal operator and
with initial point generated from the SVD of the channel matrix $C$.
The mean squared error and the averaged bit error rate are in 
Figure~\ref{fig:precoder_eq}.

\subsection{Sparse dictionary learning}
\paragraph{Problem description.}
The aim is to find a dictionary $D\in\reals^{m\times n}$ 
under which the data matrix $X\in\reals^{m\times T}$ 
can be approximated by sparse coefficients, 
i.e., $X \approx DY$ 
where $Y\in\reals^{n\times T}$ is a sparse matrix \cite{mairal2009online}.

The optimization problem can be formulated as
\[
\begin{array}{ll}
\mbox{minimize} & \frac{1}{2} \|DY-X\|_F^2 + \alpha \|Y\|_1  \\
\mbox{subject to} & \|D\|_F\leq 1,
\end{array}
\]
where the variables are $Y$ and $D$, and $\alpha>0$ is a parameter.
The problem is biconvex.

\paragraph{DMCP specification.}
The code can be written as follows.
\begin{quote}
	\begin{verbatim}
	D = Variable(m,n)
	Y = Variable(n,T)
	alpha = Parameter(sign = 'Positive')
	cost = square(norm(D*Y-X,'fro'))/2+alpha*norm(Y,1)
	prob = Problem(Minimize(cost), [norm(D,'fro') <= 1])
	prob.solve(method = 'bcd')
	\end{verbatim}
\end{quote}
\begin{figure}
	\centering
	\includegraphics[width=0.5\textwidth]{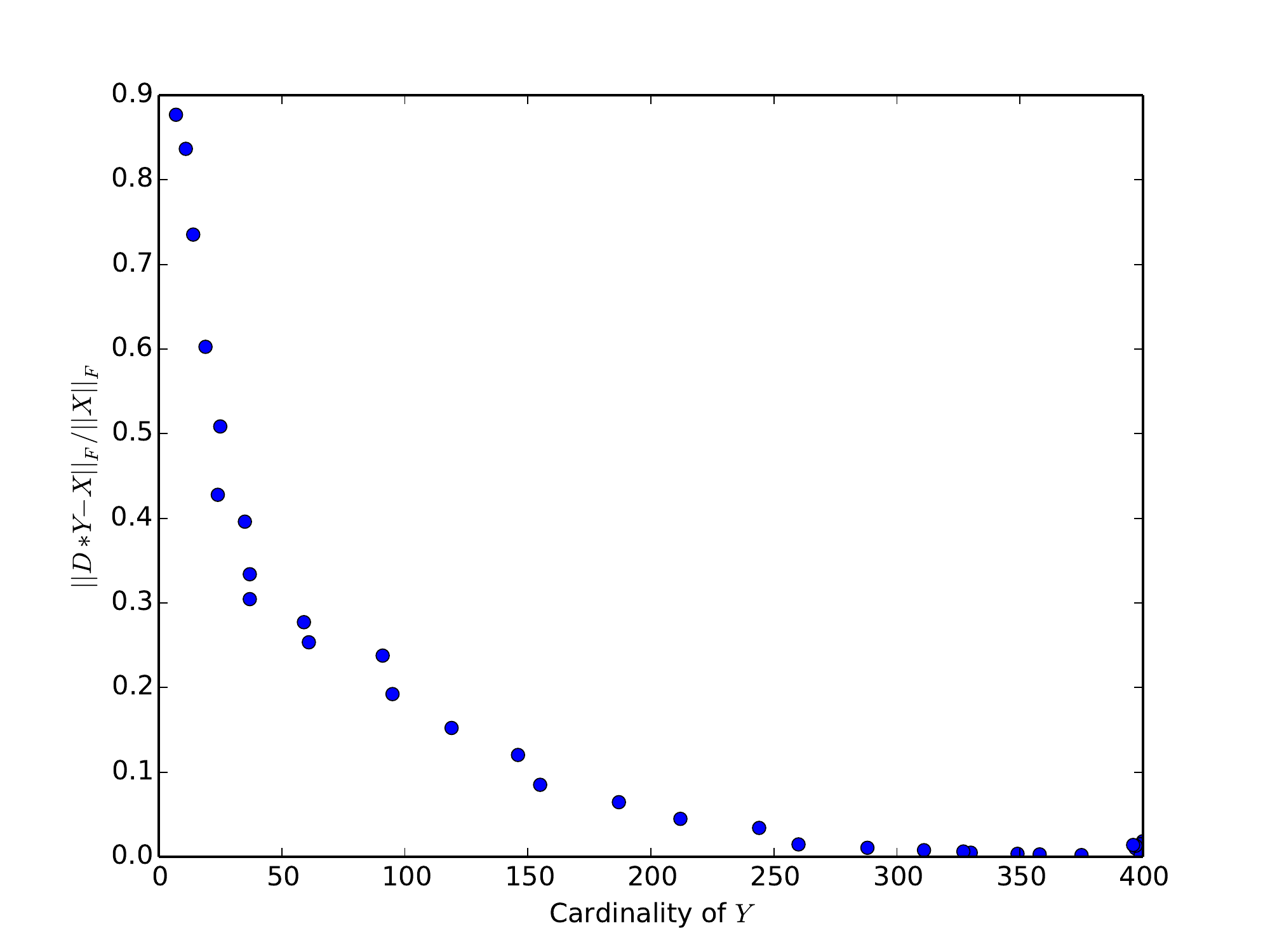}
	\caption{Sparse dictionary learning. }
	\label{fig:sparse_DL}
\end{figure}

\paragraph{Numerical result.}
In an experiment, $X$ is a random normal matrix with $m=10$, $n=20$, and $T=20$.
The parameter $\alpha$ is swept from $10^{-5}$ to $1$. For each value of 
$\alpha$, the method is called with random initialization, 
and the relative approximation error and the cardinality of $Y$ are shown 
as a blue dot in Figure \ref{fig:sparse_DL}.

\subsection{Sparse feedback matrix design}
\paragraph{Problem description.}
To design a sparse linear constant output feedback control $u = Ky$
for the system 
$$\dot{x} = Ax+Bu, \quad y = Cx,$$ 
which results in a decay rate $r$ no less than a given threshold 
$\theta>0$ in the closed-loop system, 
we consider the following optimization problem 
\cite{hassibi1999path, boyd1994linear}
\[
\begin{array}{ll}
\mbox{minimize} & \sum_{ij} |K_{ij}|  \\
\mbox{subject to} & P \succeq I,\quad  r \geq \theta \\
& -2r P \succeq (A+BKC)^T P + P(A+BKC),  
\end{array}
\]
where $K$, $P$, and $r$ are variables, 
and $A$, $B$, $C$, and $\theta$ are given.
The notation $P \succeq I$ means that $P-I$ is semidefinite.
The problem is biconvex with minimal sets of variables to fix $\{P\}$ 
and $\{K, r\}$.

\paragraph{DMCP specification.}
The code can be the following.
\begin{quote}
	\begin{verbatim}
	P = Variable(n,n)
	K = Variable(m1,m2)
	r = Variable(1)
	cost = norm(K,1)
	constr = [np.eye(n) << P, r >= theta]
	constr += [(A+B*K*C).T*P+P*(A+B*K*C) << -P*r*2]
	prob = Problem(Minimize(cost), constr)
	prob.solve(method = 'bcd')
	\end{verbatim}
\end{quote}
\paragraph{Numerical result.}
An example with $n=m_1=5$, $m_2=4$, $\theta = 0.01$, and
the following data matrices from \cite{hassibi1999path} is tested.
\begin{align*}
A &= \begin{bmatrix}
-2.45 & -0.90 & 1.53 & -1.26 & 1.76 \\
-0.12 & -0.44 & -0.01 & 0.69 & 0.90\\
2.07 & -1.20 & -1.14 & 2.04 & -0.76\\
-0.59 & 0.07 & 2.91 & -4.63 & -1.15\\
-0.74 & -0.23 & -1.19 & -0.06 & -2.52
\end{bmatrix} \\
B & = \begin{bmatrix}
0.81 & -0.79 & 0 & 0 & -0.95\\
-0.34 & -0.50 & 0.06 & 0.22 & 0.92\\
-1.32 & 1.55 & -1.22 & -0.77 & -1.14\\
-2.11 & 0.32 & 0 & -0.83 & 0.59\\
0.31 & -0.19 & -1.09 & 0 & 0
\end{bmatrix} \\
C & = \begin{bmatrix}
0 & 0 & 0.16 & 0 & -1.78 \\
1.23 & -0.38 & 0.75 & -0.38 & 0\\
0.46 & 0 & -0.05 & 0 & 0\\
0 & -0.12 & 0.23 & -0.12 & 1.14\\
\end{bmatrix} \\
\end{align*}
The initial value $P^0$ is an identity matrix, $r^0 = 1$,
and $K^0$ is an matrix with all zeros.
The result is that $r = 0.01$ and 
\begin{align*}
K = \begin{bmatrix}
0 &  0.32 &  0 & 0 \\
0 &  -0.46 &  0 &  0 \\
0 & 0 &  0 & 0 \\
0 & 0 & 0 & 0 \\
0 &  0.11 & 0 & 0 \\
\end{bmatrix},
\end{align*}
which is sparse. 
The three nonzero entries are in the second column,
so only the second output needs to be fed back.
In the work \cite{hassibi1999path} with decay rate no less than
$0.35$ another sparse feedback matrix is found
with the same cardinality $3$.

\subsection{Bilinear control}
\paragraph{Problem description.}
A discrete time $m$-input bilinear control system is of the following form
\cite{hours2014parametric}
\begin{align*}
x_{t+1} = A x_t + \sum_{i=1}^m u^i_t B^i x_t, \quad t = 1,\ldots,n-1,
\end{align*}
where $u_t = [u^1_t, \ldots, u^m_t]\in\Omega\subseteq\reals^m$ 
is the input, 
and $x_t\in\reals^{d}$ is the system state at time $t$.
In an optimal control problem with fixed initial state, 
given system matrices $A, B^i \in\reals^{n\times n}$ 
and convex objective functions $f$ and $g$,
an optimization problem can be formulated as
\[
\begin{array}{ll}
\mbox{minimize} &f(x) + g(u)\\
\mbox{subject to} 
& x_1 = x_\mathrm{ini} \\
& (x_t, u_t) \in \Omega, \quad  t = 1,\ldots,n-1\\
& x_{t+1} = A x_t + \sum_{i=1}^m u^i_t B^i x_t, \quad t = 1,\ldots,n-1,
\end{array}
\]
where $x_t$ and $u_t$ are variables,
and $\Omega$ is a given convex set describing bounds on $u_t$ and $x_t$. 
The problem is multi-convex.

As a special case,
a standard model of D.C.-motor is a bilinear system of the following form \cite{Derese1982}
\[
\dot{x} = A_0 x + u A_1 x + bv,
\]
where the derivative is with respect to time, and
\[
x = \begin{bmatrix}
x^1 \\
x^2
\end{bmatrix},\quad
A_0 = \begin{bmatrix}
-1 & 0 \\
0 & -0.1
\end{bmatrix},\quad
A_1 = \begin{bmatrix}
0 & -19 \\
0.1 & 0
\end{bmatrix},\quad
b = \begin{bmatrix}
20 \\
0
\end{bmatrix}.
\]
The correspondence between model and physical variables is that,
$x^1$ is the armature current, $x^2$ is the speed
of rotation, $u$ is the field current, and $v$ is the armature voltage.
For nominal operation,
\[
x^1 = x^2 = u = v = 1.
\]
A control problem is the braking with short-circuited armature ($v_t =0$).
The field current $u$ is controlled 
such that the rotation speed decreases to zero as fast as possible, 
and that the armature current is not excessively large. 
By discretizing over time and taking $10$ samples per second, 
the problem can be formulated in the following form
\[
\begin{array}{ll}
\mbox{minimize} & \|x^2\|_2 \\
\mbox{subject to}
&  x_1 = \begin{bmatrix} 1 \\1\end{bmatrix} \\
& \max_{t=1,\ldots, n} |x^1_t| \leq M \\
& (x_{t+1} - x_t)/0.1 = A_0 x_t + u_t A_1 x_t, 
\quad t=1,\ldots,n-1,
\end{array}
\]
where $u\in\reals^{n-1}$ and $x_t = [x^1_t, x^2_t]\in\reals^2$ 
for $t = 1,\ldots,n$ are variables,
and the notation $x^i = [x^i_1, \ldots, x^i_n]$, $i = 1,2$.
The problem is biconvex if we consider 
$x = [x_1,\ldots,x_n]\in\reals^{2\times n}$ as one variable.

\paragraph{DMCP specification.}
The code is as the following.
\begin{quote}
\begin{verbatim}
x = Variable(2,n)
u = Variable(n-1)
constr = [x[:,0] == 1, max_entries(abs(x[0,:])) <= M]
for t in range(n-1):
    constr += [x[:,t+1]-x[:,t] == 0.1*(A0*x[:,t]+A1*x[:,t]*u[t])]
prob = Problem(Minimize(norm(x[1,:])), constr)
prob.solve(method = 'bcd')
\end{verbatim}
\end{quote}

\begin{figure}
	\centering
	\includegraphics[width=0.5\textwidth]{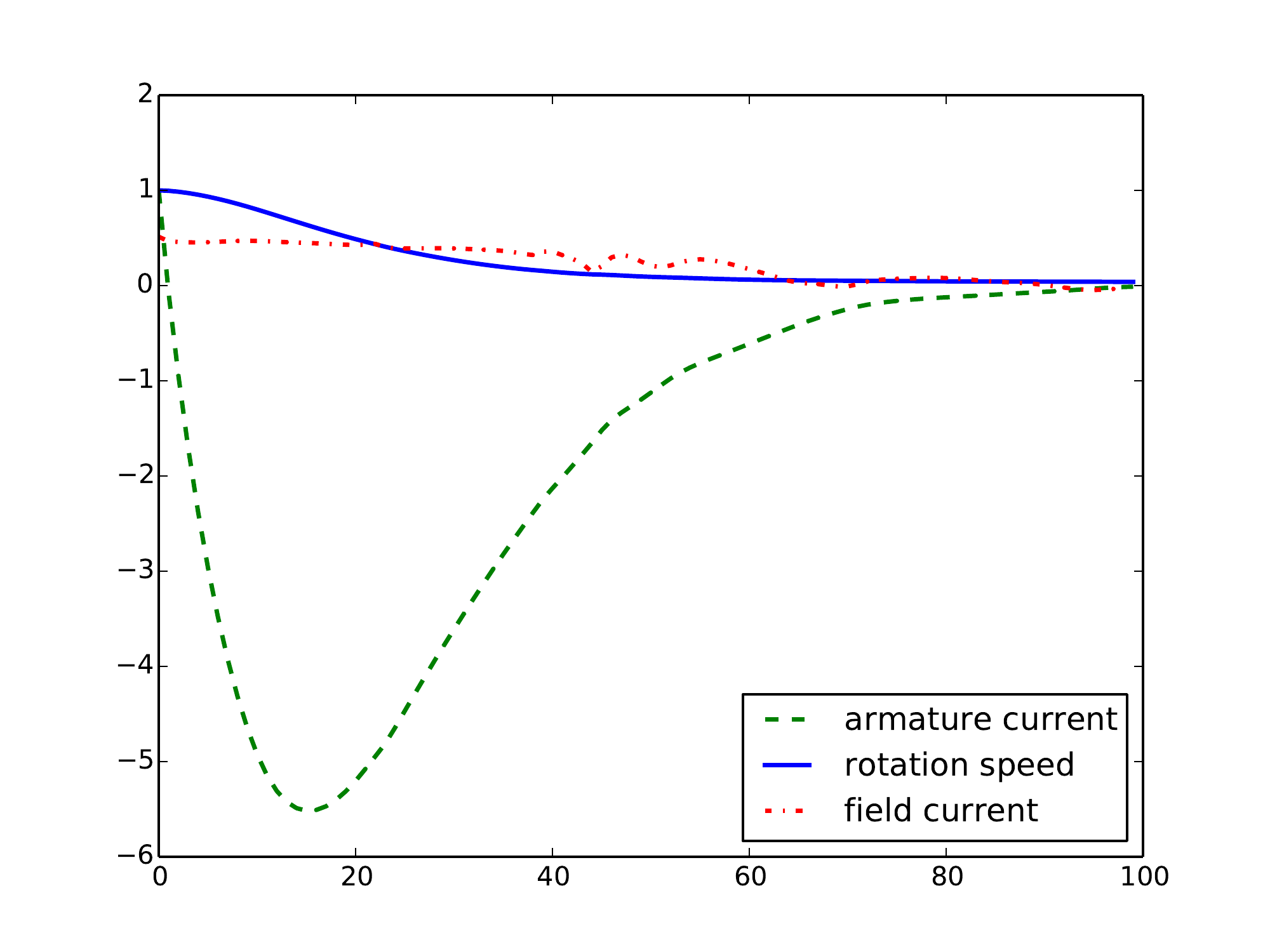}
	\caption{D.C.-motor braking.}
	\label{fig:DC_motor}
\end{figure}

\paragraph{Numerical result.}
We take an example with $n=100$ and $M=8$.
The initial value of $x$ is zero 
and of $u$ is a vector linearly decreasing from $0.5$ to $0$.
The result is shown in Fig.~\ref{fig:DC_motor},
where the braking is faster than that in a linear control system shown in
\cite{Derese1982}.

\subsection{Resistance estimation}
\paragraph{Problem description.}
A problem in direct current (DC) circuit is to estimate the values of resistors such that
certain constraints on currents and voltages can be satisfied.
The topology of the circuit is given, and several observations
on currents and voltages are known.
A general problem of estimating the resistance to fit the topology and the observations
 can be written as the following
\[
\begin{array}{ll}
\mbox{minimize} & f(u)+g(i) \\
\mbox{subject to} & A(r) i = u \\
& u \in\mathcal{U},\quad  i \in \mathcal{I},\quad r\in\mathcal{R},
\end{array}
\]
where $u\in\reals^n$, $i\in\reals^m$, and $r\in\reals_+^d$ are variables
representing voltages, currents, and resistance, respectively.
The convex functions $f$ and $g$ penalize deviations from the observations.
The first constraint corresponds to the Ohm's law, 
where the mapping $A$ is linear and depends on the topology of the circuit.
The sets $\mathcal{U}$, $\mathcal{I}$, and $\mathcal{R}$ are convex,
and they may describe the Kirchhoff's circuit laws.
The problem is multi-convex due to the first constraint.

A simple example is shown in the following circuit diagram.
\begin{center}
\begin{tikzpicture}[scale=1.4]
\draw[color=black]
(-0.2,0)node[]{$u_0$}
(0,0) to [R, l=$a_1$, -*] (0,2)
(1.5,0) to [R, l=$a_2$, *-*] (1.5,2)
(3,0) to [R, l=$a_3$, *-*] (3, 2) 
(3,0) to (3.2,0)
(3.75,0) node[]{\large{$\cdots$}}
(4.3,0) to (4.5,0) to [R, l=$a_{n-1}$, *-*] (4.5, 2)
(6,0) to [R, l=$a_n$, *-*] (6, 2)
(7.5,0) to[american current source, l=$I_0$] (7.5,4)

(0,2) to [R, l=$b_1$, *-] (0,4)
(1.5,2) to [R, l=$b_2$, *-*] (1.5,4)
(3,2) to [R, l=$b_3$, *-*] (3,4) to (3.2,4)
(3.75,4) node[]{\large{$\cdots$}}
(4.5,2) to [R, l=$b_{n-1}$, *-*] (4.5,4) to (4.3,4)
(6,2) to [R, l=$b_n$, *-*] (6,4)

(0,2) to [R, l=$c_1$](1.5,2)
(1.5,2) to [R, l=$c_2$,*-*] (3,2) to (3.2,2)
(3.75,2) node[]{\large{$\cdots$}}
(4.3,2) to (4.5,2) to [R, l=$c_{n-1}$,*-*] (6,2)

(-0.2,2.1)node[]{$v_1$}
(1.3,2.1)node[]{$v_2$}
(2.8,2.1)node[]{$v_3$}
(4.2,2.1)node[]{$v_{n-1}$}
(5.8,2.1)node[]{$v_n$}

(0.75,0.2){} node[]{$i_1$}
(2.25,0.2){} node[]{$i_2$}
(5.25,0.2){} node[]{$i_{n-1}$}
(6.75,0.2){} node[]{$i_n$}

(0.75,4.2){} node[]{$j_1$}
(2.25,4.2){} node[]{$j_2$}
(5.25,4.2){} node[]{$j_{n-1}$}
(6.75,4.2){} node[]{$j_n$}

(0,4) to node[currarrow] {} (1.5,4)
(1.5,4) to node[currarrow] {} (3,4)
(4.5,4) to node[currarrow] {} (6,4)
(6,4) to node[currarrow] {} (7.5,4)	
(7.5,4) node[ground, rotate = 90]{}

(1.5,0) to node[currarrow,rotate=180] {} (0,0)
(3,0) to node[currarrow,rotate=180] {} (1.5,0)
(6,0) to node[currarrow,rotate=180] {} (4.5,0)
(7.5,0) to node[currarrow,rotate=180] {} (6,0);	
\end{tikzpicture}
\end{center}
The known quantities are the
current source $I_0$ and the voltage level $u_0$.
It is observed that $v_k-v_{k+1} \approx \delta$ for $k = 1,\ldots,n-1$,
so the optimization problem is
\[
\begin{array}{ll}
\mbox{minimize} & \sum_{k=1}^{n-1}(v_k-v_{k+1} - \delta)^2\\
\mbox{subject to} 
&  x_{1} = y_{1}+z_{1}, \quad x_{n} + z_{n-1} = y_{n}\\
& i_1 = x_1, \quad i_n = -I_0, \quad j_1 = y_1,\quad j_n = -I_0\\
&  x_{k+1} + z_k = y_{k+1}+z_{k+1}, \quad k=1,\ldots,n-2 \\
& x_k a_k = u_0 - v_k, \quad y_k b_k = v_k, \quad k=1,\ldots,n \\
& i_{k+1} = i_k + x_{k+1}, \quad j_{k+1} = j_k + y_{k+1}, 
\quad z_k c_k = v_k - v_{k+1}, \quad k=1,\ldots,n-1,\\
\end{array}
\]
where $x, y, i, j\in\reals^n, z\in\reals^{n-1}$ are variables for currents, 
$a, b\in\reals_+^n, c\in\reals_+^{n-1}$ are the variables for resistance,
and $v\in\reals^n$ is the variable for voltages.
The problem is multi-convex, and the minimal sets to fix are not obvious.

\paragraph{DMCP specification.}
The code can be as the following.
\begin{quote}
\begin{verbatim}
x = Variable(n)
y = Variable(n)
z = Variable(n-1)
i = Variable(n)
j = Variable(n)
a = Variable(n, sign = 'Positive')
b = Variable(n, sign = 'Positive')
c = Variable(n-1, sign = 'Positive')
v = Variable(n)
constr = [x[0] == y[0]+z[0], x[n-1]+z[n-2] == y[n-1]]
constr += [i[0] == x[0], j[0] == y[0], i[n-1] == -I0, j[n-1] == -I0]
cost = 0
for k in range(n-2):
    constr += [x[k+1]+z[k] == y[k+1]+z[k+1]]
for k in range(n):
    constr += [x[k]*a[k] == u0 - v[k], y[k]*b[k] == v[k]]
for k in range(n-1):
    cost += square(v[k]-v[k+1]-delta)
    constr += [z[k]*c[k] == v[k]-v[k+1]]
    constr += [i[k+1]== i[k]+x[k+1], j[k+1] == j[k]+y[k+1]]
prob = Problem(Minimize(cost), constr)
prob.solve(method = 'bcd')
\end{verbatim}
\end{quote}
The \verb|find_minimal_sets| function returns
\begin{quote}
	\begin{verbatim}
[[5, 7, 8], [3, 5, 6], [1, 3, 7], [1, 3, 6], 
[1, 7, 8], [5, 6, 8], [1, 6, 8], [3, 5, 7]]
	\end{verbatim}
\end{quote}
where indices $1,3,5,6,7,8$ correspond to variables
$z,a,c,b,y,x$ respectively.

\paragraph{Numerical result.}
We set $n = 10, I_0 = -100, \delta = 1, u_0 = 12$,
and all variables are set with initial value $1$.
The method finds a feasible point with objective value $0$
which solves the problem globally, 
and the solution is shown in the following table.

\begin{center}
\begin{tabular}{cccc}
	$a_k$ & $b_k$ & $c_k$ & $v_k-v_{k+1}$ \\
0.120 & 1.123 & 1.800& 1.001 \\
0.224 & 0.987 & 1.776 & 1.002\\
0.329 & 0.885 & 2.249 & 1.003\\
0.431 & 0.785 & 2.904 & 1.004\\
0.535 & 0.687 & 3.653 & 1.004\\
0.641 & 0.591 & 4.066 & 1.004\\
0.743& 0.492 & 3.297 & 1.003\\
0.837 & 0.389 & 2.031 & 1.002\\
0.904 & 0.277 & 1.283 & 1.001\\
0.931& 0.149 &- &-
\end{tabular}
\end{center}

%

\subsection{Steady state of Markov chain}
\paragraph{Problem description.}
Suppose that $P_1, \ldots, P_n\in\reals^{m\times m}$ 
are $n$ transition matrices of Markov chains,
then it is known that any convex combination $P=\sum_{i=1}^n \theta_i P_i$ for 
$\theta_i \geq 0$ and $\sum_{i=1}^n\theta_i = 1$ is also a transition matrix of
a Markov chain.
Given $P_i$ for $i=1,\ldots,n$ and a convex function $f: \reals^m\to \reals$, 
the problem is to find such a convex combination, 
so that the Markov chain with respect to $P$ has a steady state vector 
$x\in\reals^m$ that achieves the minimum of $f$.

The problem can be formulated as
\[
\begin{array}{ll}
\mbox{minimize} & f(x) \\
\mbox{subject to} & P=\sum_{i=1}^n \theta_i P_i \\
& \theta_i \geq 0, \quad \sum_{i=1}^n\theta_i = 1 \\
& x_j \geq 0, \quad \sum_{j=1}^m x_j = 1 \\
& P^T x =  x,\\
\end{array}
\]
where $x$, $P$, and $\theta$ are variables.
The problem is biconvex.

As an example, $f(x) = \|x-x_0\|_2$, so the goal is 
to generate a transition matrix so that 
the steady state vector is close to a given distribution $x_0$.

\paragraph{DMCP specification.}
The code is as follows.
\begin{quote}
\begin{verbatim}
cost = norm(x-x0)
constr = [theta >= 0, sum_entries(theta) == 1, x >= 0, sum_entries(x) == 1]
right = 0
for i in range(n):
    right += theta[i]*P0[i]
constr += [P == right, P.T*x == x]
prob = Problem(Minimize(cost), constr)
prob.solve(method = 'bcd')
\end{verbatim}
\end{quote}

\paragraph{Numerical result.}
An example with $n=4$, $m=3$, randomly generated $P_i$, and
\[
x_0 = \begin{bmatrix}
0.25\\
0.3\\
0.45
\end{bmatrix},
\]
is tested.
The initial values are random.
The result gives
\[
x= \begin{bmatrix}
0.25\\
0.3\\
0.45
\end{bmatrix},\quad
\theta = \begin{bmatrix}
0.27\\
0.61\\
0.04\\
0.08
\end{bmatrix},
\]
which achieves the targeted steady state vector.

\subsection{Blind deconvolution}
\paragraph{Problem description.}
Blind deconvolution is an inverse problem commonly encountered in many 
practical applications such as image restoration, system identification, 
and channel estimation.
The problem is to find two vectors with some priors, 
such that their convolution approximates the given data.

Suppose that the data $d\in\reals^{m+n-1}$ is the convolution of 
an unknown sparse vector $x_0\in\reals^n$ and an unknown vector $y_0\in\reals^m$.
A problem can be formulated as the following
\[
\begin{array}{ll}
\mbox{minimize} & \|x * y-d\|_2+\alpha \|x\|_1\\
\mbox{subject to} & \|y\|_\infty \leq M,
\end{array}
\]
where $x\in\reals^n$ and $y\in\reals^m$ are variables,
and $\alpha>0$ is a parameter.
The problem is biconvex.
For any $x*y = d$ and a positive scalar $k$
the convolution of $kx$ and $y/k$ is also $d$,
so the constraint $\|y\|_\infty \leq M$ is needed to exclude a trivial solution
$x\approx 0$.

\paragraph{DMCP specification.}
The code can be written as the following.
\begin{quote}
	\begin{verbatim}
	y = Variable(m)
	x = Variable(n)
	cost = norm(conv(y,x)-d,2) + alpha*norm(x,1)
	prob = Problem(Minimize(cost), [norm(y,'inf') <= M])
	prob.solve(method = 'bcd')
	\end{verbatim}
\end{quote}

\begin{figure}
	\centering
		\centering
		\includegraphics[height=1.8in]{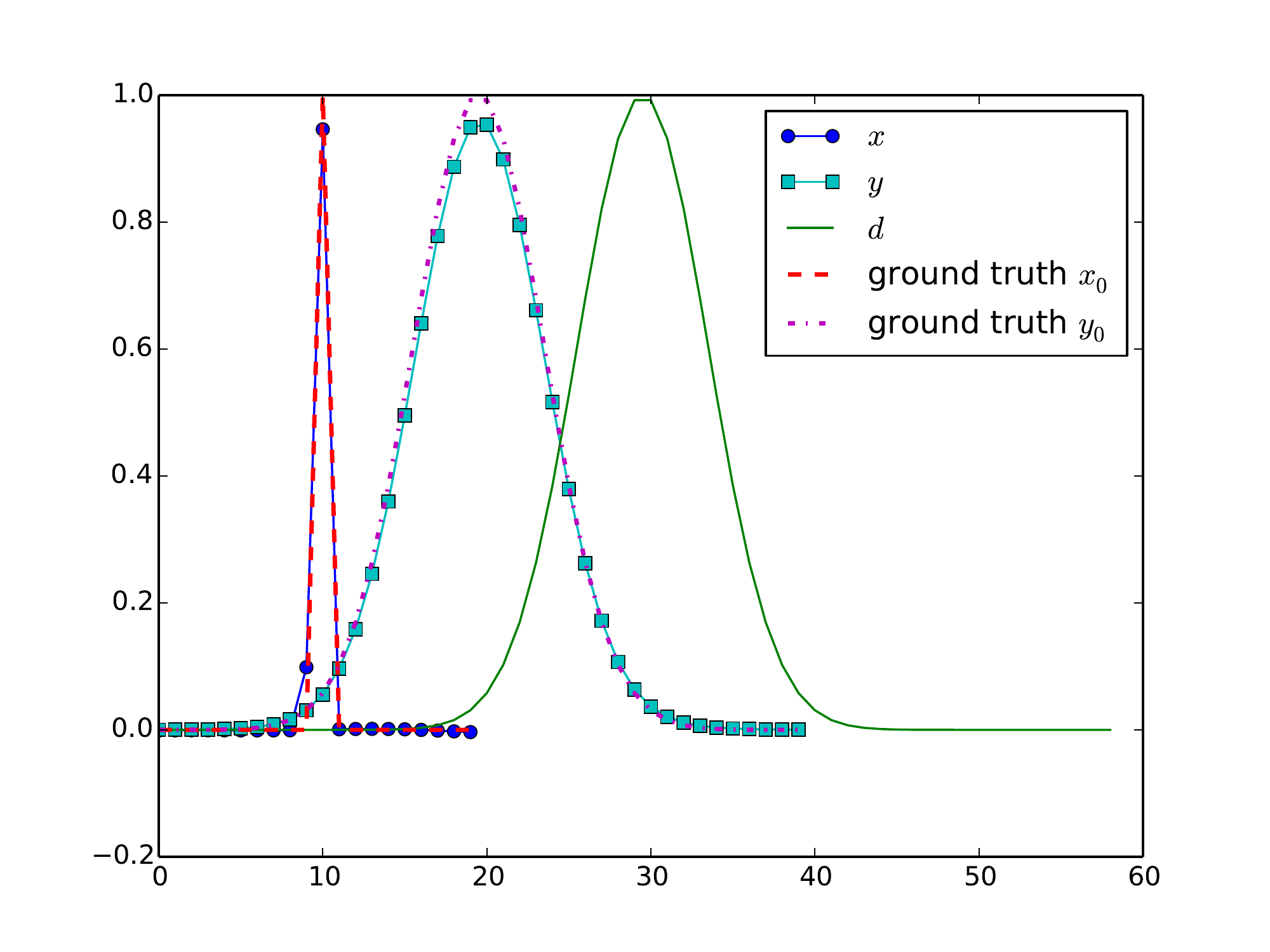}
	\caption{Blind deconvolution.}
	\label{fig:blind_dconv}
\end{figure}

\paragraph{Numerical result.}
In an example, $m=100$, $n=40$, $M=10$, $\alpha = 0.28$,
and the initial value of every variable is a vector of all ones.
The result is in Figure~\ref{fig:blind_dconv}.

\subsection*{Acknowledgments}
This material is based
upon work supported by the National Science Foundation Graduate Research Fellowship
under Grant No. DGE-114747, by the DARPA X-DATA and SIMPLEX programs,
and by the CSC State Scholarship  Fund.
\clearpage

\bibliography{dmcp}

\newcommand{\etalchar}[1]{$^{#1}$}
\begin{thebibliography}{BEGFB94}

\bibitem[ABRS10]{attouch2010}
H.~Attouch, J.~Bolte, P.~Redont, and A.~Soubeyran.
\newblock Proximal alternating minimization and projection methods for
  nonconvex problems: An approach based on the {K}urdyka-{L}ojasiewicz
  inequality.
\newblock {\em Mathematics of Operations Research}, 35(2):438--457, 2010.

\bibitem[BEGFB94]{boyd1994linear}
S.~Boyd, L.~El~Ghaoui, E.~Feron, and V.~Balakrishnan.
\newblock {\em Linear matrix inequalities in system and control theory},
  volume~15.
\newblock SIAM, 1994.

\bibitem[BV04]{BoV:04}
S.~Boyd and L.~Vandenberghe.
\newblock {\em Convex Optimization}.
\newblock Cambridge University Press, 2004.

\bibitem[CDPR04]{Chu04}
M.~Chu, F.~Diele, R.~Plemmons, and S.~Ragni.
\newblock Optimality, computation, and interpretation of nonnegative matrix
  factorizations.
\newblock {\em SIAM JOURNAL ON MATRIX ANALYSIS}, pages 4--8030, 2004.

\bibitem[DB16]{cvxpy_paper}
S.~Diamond and S.~Boyd.
\newblock {CVXPY}: A {P}ython-embedded modeling language for convex
  optimization.
\newblock {\em To appear, Journal of Machine Learning Research}, 2016.

\bibitem[DN82]{Derese1982}
I.~Derese and E.~Noldus.
\newblock Stabilization of bilinear systems.
\newblock In {\em Analysis and Optimization of Systems: Proceedings of the
  Fifth International Conference on Analysis and Optimization of Systems
  Versailles, December 14--17, 1982}, pages 974--987, 1982.

\bibitem[GB08]{GB:08}
M.~Grant and S.~Boyd.
\newblock Graph implementations for nonsmooth convex programs.
\newblock In V.~Blondel, S.~Boyd, and H.~Kimura, editors, {\em Recent Advances
  in Learning and Control}, Lecture Notes in Control and Information Sciences,
  pages 95--110. Springer, 2008.

\bibitem[GB14]{cvx}
M.~Grant and S.~Boyd.
\newblock {CVX}: {MATLAB} software for disciplined convex programming, version
  2.1.
\newblock \url{http://cvxr.com/cvx}, March 2014.

\bibitem[GBY06]{GBY:06}
M.~Grant, S.~Boyd, and Y.~Ye.
\newblock Disciplined convex programming.
\newblock In L.~Liberti and N.~Maculan, editors, {\em Global Optimization: From
  Theory to Implementation}, Nonconvex Optimization and its Applications, pages
  155--210. Springer, 2006.

\bibitem[HHB99]{hassibi1999path}
A.~Hassibi, J.~How, and S.~Boyd.
\newblock A path-following method for solving {BMI} problems in control.
\newblock In {\em American Control Conference, 1999. Proceedings of the 1999},
  volume~2, pages 1385--1389. IEEE, 1999.

\bibitem[HJ14]{hours2014parametric}
J.~Hours and C.~Jones.
\newblock A parametric multiconvex splitting technique with application to
  real-time {NMPC}.
\newblock In {\em 53rd IEEE Conference on Decision and Control}, pages
  5052--5057. IEEE, 2014.

\bibitem[KP08]{Kim2008}
H.~Kim and H.~Park.
\newblock Nonnegative matrix factorization based on alternating nonnegativity
  constrained least squares and active set method.
\newblock {\em SIAM Journal on Matrix Analysis and Applications},
  30(2):713--730, 2008.

\bibitem[Lof04]{Lofberg:04}
J.~Lofberg.
\newblock {YALMIP}: A toolbox for modeling and optimization in {MATLAB}.
\newblock In {\em Proceedings of the {IEEE} International Symposium on Computed
  Aided Control Systems Design}, pages 294--289, September 2004.

\bibitem[LS99]{lee1999learning}
D.~D. Lee and H.~S. Seung.
\newblock Learning the parts of objects by non-negative matrix factorization.
\newblock {\em Nature}, 401(6755):788--791, 1999.

\bibitem[LS00]{Lee00}
D.~D. Lee and H.~S. Seung.
\newblock Algorithms for non-negative matrix factorization.
\newblock In {\em In NIPS}, pages 556--562. MIT Press, 2000.

\bibitem[MBPS09]{mairal2009online}
J.~Mairal, F.~Bach, J.~Ponce, and G.~Sapiro.
\newblock Online dictionary learning for sparse coding.
\newblock In {\em Proceedings of the 26th annual international conference on
  machine learning}, pages 689--696. ACM, 2009.

\bibitem[NN92]{NesNem:92}
Y.~Nesterov and A.~Nemirovsky.
\newblock Conic formulation of a convex programming problem and duality.
\newblock {\em Optimization Methods and Software}, 1(2):95--115, 1992.

\bibitem[NW06]{nocedal2006numerical}
J.~Nocedal and S.~Wright.
\newblock {\em Numerical optimization}.
\newblock Springer Science \& Business Media, 2006.

\bibitem[PB14]{parikh2014proximal}
N.~Parikh and S.~Boyd.
\newblock Proximal algorithms.
\newblock {\em Foundations and Trends in optimization}, 1(3):127--239, 2014.

\bibitem[Pow73]{Powell1973}
M.~J.~D. Powell.
\newblock On search directions for minimization algorithms.
\newblock {\em Mathematical Programming}, 4(1):193--201, 1973.

\bibitem[RHL12]{razaviyayn2012unified}
M.~Razaviyayn, M.~Hong, and Z.~Luo.
\newblock A unified convergence analysis of coordinatewise successive
  minimization methods for nonsmooth optimization.
\newblock {\em preprint}, 2012.

\bibitem[Roc15]{rockafellar}
Ralph~Tyrell Rockafellar.
\newblock {\em Convex analysis}.
\newblock Princeton university press, 2015.

\bibitem[SGL94]{751690}
M.~G. Safonov, K.~C. Goh, and J.~H. Ly.
\newblock Control system synthesis via bilinear matrix inequalities.
\newblock In {\em American Control Conference, 1994}, volume~1, pages 45--49
  vol.1, June 1994.

\bibitem[SY04]{1254038}
S.~Serbetli and A.~Yener.
\newblock Transceiver optimization for multiuser {MIMO} systems.
\newblock {\em IEEE Transactions on Signal Processing}, 52(1):214--226, Jan
  2004.

\bibitem[THB15]{7170815}
J.~Thai, R.~Hariss, and A.~Bayen.
\newblock A multi-convex approach to latency inference and control in traffic
  equilibria from sparse data.
\newblock In {\em 2015 American Control Conference (ACC)}, pages 689--695, July
  2015.

\bibitem[Tse93]{Tseng1993}
P.~Tseng.
\newblock Dual coordinate ascent methods for non-strictly convex minimization.
\newblock {\em Mathematical Programming}, 59(1):231--247, 1993.

\bibitem[Tse01]{Tseng2001}
P.~Tseng.
\newblock Convergence of a block coordinate descent method for
  nondifferentiable minimization.
\newblock {\em Journal of Optimization Theory and Applications},
  109(3):475--494, 2001.

\bibitem[UHZB14]{udell2014generalized}
M.~Udell, C.~Horn, R.~Zadeh, and S.~Boyd.
\newblock Generalized low rank models.
\newblock {\em arXiv preprint arXiv:1410.0342}, 2014.

\bibitem[UMZ{\etalchar{+}}14]{cvxjl}
M.~Udell, K.~Mohan, D.~Zeng, J.~Hong, S.~Diamond, and S.~Boyd.
\newblock Convex optimization in {J}ulia.
\newblock In {\em Proceedings of the Workshop for High Performance Technical
  Computing in Dynamic Languages}, pages 18--28, 2014.

\bibitem[War63]{Warga}
J.~Warga.
\newblock Minimizing certain convex functions.
\newblock {\em Journal of the Society for Industrial and Applied Mathematics},
  11(3):588--593, 1963.

\bibitem[WYZ12]{Wen2012}
Z.~Wen, W.~Yin, and Y.~Zhang.
\newblock Solving a low-rank factorization model for matrix completion by a
  nonlinear successive over-relaxation algorithm.
\newblock {\em Mathematical Programming Computation}, 4(4):333--361, 2012.

\bibitem[XY13]{xu2013block}
Y.~Xu and W.~Yin.
\newblock A block coordinate descent method for regularized multiconvex
  optimization with applications to nonnegative tensor factorization and
  completion.
\newblock {\em SIAM Journal on imaging sciences}, 6(3):1758--1789, 2013.

\end{thebibliography}
\end{document}